\documentclass[11pt,oneside,a4paper]{article}
\usepackage{blindtext}
\usepackage{amsmath,amsfonts,amssymb,amsthm,eucal}
\usepackage{microtype}
\usepackage[colorlinks, citecolor=blue, linkcolor=blue, urlcolor=Maroon, filecolor=Maroon, linktocpage=true]{hyperref}
\usepackage[usenames,dvipsnames,svgnames,table]{xcolor}
\usepackage[normalem]{ulem}
\usepackage{lmodern}
\usepackage[inline]{enumitem}
\usepackage{bbm}
\usepackage{tikz}
\usepackage{tikz-cd}
\usepackage[font=footnotesize,labelfont=bf]{caption}
\usepackage{subcaption}
\usepackage{framed}
\usepackage[nottoc,notlot,notlof]{tocbibind}
\usepackage{tocloft}
\usepackage{morefloats}
\usepackage{float}
\usepackage[nottoc,notlot,notlof]{tocbibind}
\usepackage{tocloft}
\usepackage[backgroundcolor=LightBlue,linecolor=LightBlue]{todonotes}
\usepackage{hhline}
\usepackage{array,multirow}
\usepackage{fancyvrb}
\usepackage{listings}
\usepackage{longtable}
\usepackage{authblk}

\newcommand{\tr}{\,\mathrm{tr}}

\theoremstyle{definition}
\newtheorem{Th}{Theorem}[section]

\newtheorem{Prop}[Th]{Proposition}

\newtheorem{Lem}[Th]{Lemma}

\theoremstyle{definition}
\newtheorem{rem}[Th]{Remark}
\usepackage[a4paper,
	hcentering,
	text={453pt,686pt},
]{geometry}
\allowdisplaybreaks
\catcode`@=11
\@addtoreset{equation}{section}
\catcode`@=12
\begin{document}
\title{\textbf{\Large Classification of left-invariant Einstein metrics on $\mathrm{SL}(2,\mathbb{R})\times\mathrm{SL}(2,\mathbb{R})$ that are bi-invariant under a one-parameter subgroup}}
\author[1]{Vicente Cort\'es}
\author[1]{Jeremias Ehlert}
\author[1]{Alexander S. Haupt}
\author[2]{David Lindemann}
\affil[1]{Department of Mathematics and Center for Mathematical Physics\linebreak
University of Hamburg\linebreak
Bundesstra{\ss}e 55, D-20146 Hamburg, Germany\linebreak
\texttt{vicente.cortes@uni-hamburg.de, jeremias.ehlert@gmx.de, a-haupt@web.de}\linebreak}
\affil[2]{Department of Mathematics\linebreak
Aarhus University\linebreak
Ny Munkegarde 118, DK-8000 Aarhus C, Denmark\linebreak
\texttt{david.lindemann@math.au.dk}}
\date{}
\maketitle
\begin{abstract}
	\noindent
	We classify all left-invariant pseudo-Riemannian Einstein metrics on $\mathrm{SL}(2,\mathbb{R})\times \mathrm{SL}(2,\mathbb{R})$ that are bi-invariant under a one-parameter subgroup. We find that there are precisely two such metrics up to homothety, the Killing form and a nearly pseudo-K\"ahler metric.
\end{abstract}
\textbf{Keywords:} special linear group, Einstein manifolds, non-compact homogeneous pseudo-Riemannian manifolds\\
\textbf{MSC classification:} 53C25, 53C30 (primary)
\tableofcontents
\section{Introduction}
	In recent work by C. B\"ohm and R.A. Lafuente \cite{BL1} it has been proven that for all integers $k\geq 2$, $\mathrm{SL}(2,\mathbb{R})^k$ does not admit a Riemannian Einstein metric, which is consistent with our main Theorem \ref{thm_main}. This result, however, has recently been fundamentally extended by C. B\"ohm and R.A. Lafuente in \cite{BL2} in which the authors prove the long-standing Alekseevsky conjecture, stating that every connected homogeneous Einstein manifold of negative scalar curvature is diffeomorphic to $\mathbb{R}^n$. Before this was shown, the Alekseevsky conjecture was known to hold up to dimension $5$, and in dimension $6$ the only open case was $\mathrm{SL}(2,\mathbb{C})$, cf. \cite{B}. If we allow for the metric to be pseudo-Riemannian, the classification of left-invariant Einstein metrics on $\mathrm{SL}(2,\mathbb{R})^k$ for $k\geq 2$ is however still an open problem. In this work we will classify all left-invariant pseudo-Riemannian Einstein metrics on $\mathrm{SL}(2,\mathbb{R})\times \mathrm{SL}(2,\mathbb{R})$ that are additionally bi-invariant under a one-parameter subgroup. We will proceed by first classifying all generators of one-parameter subgroups up to automorphisms and the corresponding invariant pseudo-Riemannian metrics, see Proposition \ref{prop_classification_inv_subgrp_02tensors}. Then, with the help of computer algebra software using Gr\"obner bases, we solve the Einstein equations in each case. The results are displayed in Lemma \ref{lem_isometry_class_solutions}, and in Lemmas \ref{lem_family_solutions_representative}, \ref{lem_equiv_B_g1}, \ref{lem_inequivalence_metrics}, and \ref{lem_widetildeg_nK} we determine the solutions up to automorphisms. We obtain the following as our main result.
	
	\begin{Th}\label{thm_main}
		Let $g$ be a left-invariant pseudo-Riemannian Einstein metric on $\mathrm{SL}(2,\mathbb{R})\times\mathrm{SL}(2,\mathbb{R})$ that is bi-invariant under a one-parameter subgroup. Then $g$ is either homothetic to the Killing form $B$, or to the up to homothety unique nearly pseudo-K\"ahler metric given in the basis \eqref{eqn_basis_sl2_sl2} at the neutral element by
			\begin{equation*}
				\widetilde{g}:=\dfrac{1}{\sqrt{3}}\left(
					\begin{array}{ccc|ccc}
						-2 & 0 & 0 & 1 & 0 & 0\\
						0 & 2 & 0 & 0 & 1 & 0\\
						0 & 0 & 2 & 0 & 0 & 1\\
						\hline
						\overset{\vphantom{.}}{1} & 0 & 0 & -2 & 0 & 0\\
						0 & 1 & 0 & 0 & 2 & 0\\
						0 & 0 & 1 & 0 & 0 & 2
					\end{array}
				\right).
			\end{equation*}
		Positive rescalings of $\widetilde{g}$ and the Killing form are pseudo-Riemannian with index $2$.
	\end{Th}
	
	The existence of a nearly pseudo-K\"ahler metric on $\mathrm{SL}(2,\mathbb{R})\times\mathrm{SL}(2,\mathbb{R})$ follow from \cite{G}. Its uniqueness up to homothety has been proven in \cite{SS}. In order to see that the left-invariant metric corresponding to $\widetilde{g}$ is nearly pseudo-K\"ahler, one checks that it is obtained by the pullback of the Killing form $B\times B\times B$ on $\mathrm{SL}(2,\mathbb{R})\times\mathrm{SL}(2,\mathbb{R})\times\mathrm{SL}(2,\mathbb{R})$ to the diagonal embedding of $\mathrm{SL}(2,\mathbb{R})$, $(\mathrm{SL}(2,\mathbb{R})\times\mathrm{SL}(2,\mathbb{R})\times\mathrm{SL}(2,\mathbb{R}))/\Delta$, which we identify with $\mathrm{SL}(2,\mathbb{R})\times\mathrm{SL}(2,\mathbb{R})$. On the level of Lie algebras, this identification corresponds to
		\begin{align*}
			\iota:\mathfrak{sl}(2,\mathbb{R})\times \mathfrak{sl}(2,\mathbb{R}) &\to T_e (\mathrm{SL}(2,\mathbb{R})\times\mathrm{SL}(2,\mathbb{R})\times\mathrm{SL}(2,\mathbb{R})/\Delta)\subset\mathfrak{sl}(2,\mathbb{R})\times\mathfrak{sl}(2,\mathbb{R})\times\mathfrak{sl}(2,\mathbb{R}),\\
			(X,Y) &\mapsto (X,Y,-X-Y),
		\end{align*}
	where $T_e (\mathrm{SL}(2,\mathbb{R})\times\mathrm{SL}(2,\mathbb{R})\times\mathrm{SL}(2,\mathbb{R})/\Delta)$ is identified with $\{(X,Y,Z)\in \mathfrak{sl}(2,\mathbb{R})\times\mathfrak{sl}(2,\mathbb{R})\times\mathfrak{sl}(2,\mathbb{R})\ |\ X+Y+Z=0\}$. When taking the same basis in each factor of $\mathfrak{sl}(2,\mathbb{R})\times \mathfrak{sl}(2,\mathbb{R})$, the pullback of $B\times B\times B$ is of the form
		\begin{equation*}
			\iota^*(B\times B\times B)=\left(
				\begin{array}{c|c}
					2B & B\\
					\hline
					\overset{\vphantom{.}}{B} & 2B
				\end{array}
			\right).
		\end{equation*}
	We see, however, that the above form does not coincide with $\widetilde{g}$. The reason is that the induced Lie algebra structure on $\mathfrak{sl}(2,\mathbb{R})\times \mathfrak{sl}(2,\mathbb{R})$ in the basis \eqref{eqn_E_basis} in each factor of the image of $\iota$ is isomorphic to, but does not coincide with the Lie algebra structure on $\mathfrak{sl}(2,\mathbb{R})\times \mathfrak{sl}(2,\mathbb{R})$ in the basis \eqref{eqn_E_basis} in each of the two factors. In fact, one can show that the structure constants $\overline{c}^k_{ij}$ on $\mathfrak{sl}(2,\mathbb{R})\times \mathfrak{sl}(2,\mathbb{R})$ of the induced Lie algebra structure in any orthonormal basis of $(\mathfrak{sl}(2,\mathbb{R})\times \mathfrak{sl}(2,\mathbb{R}),B\times B)$ that is obtained by taking a copy of an orthonormal basis of $(\mathfrak{sl}(2,\mathbb{R}),B)$ in each factor are of the form
		\begin{align*}
			&(\overline{c}^k_{ij})_{ij}=\left(
				\begin{array}{c|c}
					\underset{}{\vphantom{.}}\frac{1}{3}c^k_{ij} & -\frac{1}{3}c^k_{ij}\\
					\hline
					\overset{}{\vphantom{\big(}}\frac{1}{3}c^k_{ij} & -\frac{2}{3}c^k_{ij}
				\end{array}
			\right),\ 1\leq k \leq 3,\\
			&(\overline{c}^k_{ij})_{ij}=\left(
				\begin{array}{c|c}
					\underset{}{\vphantom{.}}-\frac{2}{3}c^k_{ij} & -\frac{1}{3}c^k_{ij}\\
					\hline
					\overset{}{\vphantom{\big(}}\frac{1}{3}c^k_{ij} & \frac{1}{3}c^k_{ij}
				\end{array}
			\right),\ 4\leq k \leq 6,
		\end{align*}
	where the $c^k_{ij}$ denote the structure constants of $\mathfrak{sl}(2,\mathbb{R})$. An analogous formula also holds if we replace $\mathfrak{sl}(2,\mathbb{R})$ by any other semi-simple Lie algebra. We will not give an explicit transformation to show that $\widetilde{g}$ is nearly pseudo-K\"ahler, but instead use the uniqueness result from \cite{SS}, see Lemma \ref{lem_widetildeg_nK}.
	
	Theorem \ref{thm_main} is particularly interesting when compared with a result by Nikonorov and Rodionov. In \cite[Thm.\,2]{NR} the authors show that there exist precisely two left-invariant Einstein metrics on $\mathrm{SU}(2)\times\mathrm{SU}(2)$ that are bi-invariant under an embedding of $S^1\subset\mathrm{SU}(2)\times\mathrm{SU}(2)$ up to homothety. One of them is the standard product metric, and the other a nearly K\"ahler metric first described in \cite{J}. Their result was generalised in \cite{BCHL} to only require bi-invariance under a discrete subgroup with at least two elements and non-trivial action on $S^1\subset\mathrm{SU}(2)\times\mathrm{SU}(2)$. Note that both $\mathrm{SU}(2)\times\mathrm{SU}(2)$ and $\mathrm{SL}(2,\mathbb{R})\times\mathrm{SL}(2,\mathbb{R})$ are real forms of $\mathrm{SL}(2,\mathbb{C})\times\mathrm{SL}(2,\mathbb{C})$. Interestingly, we find in Theorem \ref{thm_main} that for $\mathrm{SL}(2,\mathbb{R})\times\mathrm{SL}(2,\mathbb{R})$ an analogous result to \cite[Thm.\,2]{NR} holds, that is that under the assumption of bi-invariance under a one-parameter subgroup there exist up to homothety also only two left-invariant Einstein metrics on $\mathrm{SL}(2,\mathbb{R})\times\mathrm{SL}(2,\mathbb{R})$, one being the standard product metric and the other the nearly K\"ahler metric. It is however still an open question whether our result has a similar generalisation to metrics invariant under discrete subgroups.
	
	In comparison with a previous work of three of the authors in \cite{BCHL} we face similar problems when attempting to solve the Einstein equations in the most general case, that is without further restricting the left-invariant pseudo-Riemannian metric. This is mainly due to the fact that our currently available computers are too slow as parallelising known algorithms for finding Gr\"obner bases is a highly non-trivial problem \cite{AGLZ,MP}. Furthermore, the computation of Gr\"obner bases tends to use a lot of RAM, cf. for example \cite{BCHL} for the details of a problem with similar complexity.
	
\paragraph*{Acknowledgements}\textcolor{white}{.}\\
	\noindent
	This work was supported by the German Science Foundation (DFG) under Germany's Excellence Strategy -- EXC 2121 ``Quantum Universe'' -- 390833306 and partly supported by the German Science Foundation (DFG) under the Research Training Group 1670 ``Mathematics inspired by String Theory''.

\section{Preliminaries}
	In the following we will give a quick overview of the notation and known results that we will be using in our proof of Theorem \ref{thm_main}. Unless explicitly stated otherwise, we use the Einstein summation convention. We will work with the following basis of $\mathfrak{sl}(2,\mathbb{R})$ in our calculations
		\begin{equation}\label{eqn_E_basis}
			E_1:=\left(\begin{array}{cc} 0 & -1\\ 1 & 0\end{array}\right),\quad E_2:=\left(\begin{array}{cc} 0 & 1\\ 1 & 0\end{array}\right),\quad E_3:=\left(\begin{array}{cc} 1 & 0\\ 0 & -1\end{array}\right),
		\end{equation}
	so that $[E_1,E_2]=-2E_3$, $[E_1,E_3]=2E_2$, and $[E_2,E_3]=2E_1$. We further define
		\begin{equation*}
			N:=\frac{1}{2}(E_2-E_1)=\left(\begin{array}{rr} 0 & 1\\ 0 & 0\end{array}\right).
		\end{equation*}
	In our basis we thus have
		\begin{equation*}
			\mathrm{ad}_{E_1}=\left(\begin{array}{ccc} 0 & 0 & 0\\ 0 & 0 & 2\\ 0 & -2 & 0\end{array}\right),\quad \mathrm{ad}_{E_2}=\left(\begin{array}{ccc} 0 & 0 & 2\\ 0 & 0 & 0\\ 2 & 0 & 0\end{array}\right),\quad \mathrm{ad}_{E_3}=\left(\begin{array}{ccc} 0 & -2 & 0\\ -2 & 0 & 0\\ 0 & 0 & 0\end{array}\right),
		\end{equation*}
	and
		\begin{equation*}
			\mathrm{ad}_{N}=\left(\begin{array}{ccc} 0 & 0 & 1\\ 0 & 0 & -1\\ 1 & 1 & 0\end{array}\right).
		\end{equation*}
	The Killing form of $\mathrm{SL}(2,\mathbb{R})$, $B(X,Y)=4\tr(XY)$, is given by $B=8\left(-(E_1^*)^2 + (E_2^*)^2 + (E_3^*)^2\right)$. For the Lie algebra of $\mathrm{SL}(2,\mathbb{R})\times \mathrm{SL}(2,\mathbb{R})$ we will work with the basis
		\begin{equation}
			F_1=(E_1,0),\ F_2=(E_2,0),\ F_3=(E_3,0),\ F_4=(0,E_1),\ F_5=(0,E_2),\ F_6=(0,E_3).\label{eqn_basis_sl2_sl2}
		\end{equation}
	We denote the corresponding structure constants of $\mathfrak{sl}(2,\mathbb{R})\oplus \mathfrak{sl}(2,\mathbb{R})$ in the above basis by $c_{ij}^k$, $1\leq i,j,k\leq 6$. In the following, we will identify all left-invariant tensor fields on $\mathrm{SL}(2,\mathbb{R})\times \mathrm{SL}(2,\mathbb{R})$ with their value at the neutral element. Let $g$ be a left-invariant pseudo-Riemannian metric on $\mathrm{SL}(2,\mathbb{R})\times \mathrm{SL}(2,\mathbb{R})$. In the basis $\{F_i\ |\ 1\leq i\leq 6\}$, the connection coefficients $\omega_{ij}^k$ and the Riemann curvature tensor $R$ at the neutral element are given by
		\begin{align}
			\omega^k_{ij}&=\frac{1}{2}\left(-g_{il}g^{mk} c^l_{jm} - g_{jl}g^{mk} c^l_{im} + c^k_{ij}\right),\label{eqn_con_coeff}\\
			R^l{}_{ijk}&=\omega^m_{jk}\omega^l_{im} - \omega^m_{ik}\omega^l_{jm} - c^m_{ij}\omega^l_{mk},\label{eqn_riemann_curvature}
		\end{align}
	where $R^l{}_{ijk}=F_l^*(R(F_i,F_j)F_k)$, $g_{ij}=g(F_i,F_j)$, and $g^{ij}=(g^{-1})(F^*_i,F^*_j)$ for all $1\leq i,j,k,l\leq 6$. The components of the Ricci tensor $\mathrm{Ric}$ at the neutral element are $\mathrm{Ric}_{ij}=R^l{}_{lij}$, and the Einstein equation with Einstein constant $\lambda\in\mathbb{R}$ is given by $\mathrm{Ric}=\lambda g$.

	Our main theorem is not the most general classification of all pseudo-Riemannian Einstein metrics on $\mathrm{SL}(2,\mathbb{R})\times \mathrm{SL}(2,\mathbb{R})$, but assumes invariance under a one-parameter subgroup. In the next step, we will classify all such one-parameter subgroups up to automorphisms of $\mathrm{SL}(2,\mathbb{R})\times \mathrm{SL}(2,\mathbb{R})$. Note that automorphisms preserve the structure constants and the Killing form, and recall that the automorphism group of $\mathrm{SL}(2,\mathbb{R})$ is the projective linear group 
		\begin{equation*}
			\mathrm{Aut}(\mathrm{SL}(2,\mathbb{R}))=\mathrm{PGL}(2,\mathbb{R}).
		\end{equation*}
	It follows that
		\begin{equation*}
			\mathrm{Aut}(\mathrm{SL}(2,\mathbb{R})\times \mathrm{SL}(2,\mathbb{R}))=\left(\mathrm{PGL}(2,\mathbb{R})\times \mathrm{PGL}(2,\mathbb{R})\right)\rtimes \mathbb{Z}_2,
		\end{equation*}
	where the $\mathbb{Z}_2$-factor corresponds to swapping the elements of the two copies of $\mathrm{SL}(2,\mathbb{R})$.
		
		\begin{Lem}\label{lem_sl2_under_auto_class}
			Under the induced action of $\mathrm{PGL}(2,\mathbb{R})$, each non-zero element in $\mathfrak{sl}(2,\mathbb{R})$ is equivalent to either $rE_1$ or $rE_3$ for precisely one $r>0$, or to $N$.
			\begin{proof}
				The Lie algebra $\mathfrak{sl}(2,\mathbb{R})$ consists of traceless $2\times 2$-matrices with real entries. Hence, if $\lambda$ is a solution of its characteristic polynomial, so is its conjugate $\overline{\lambda}$. The traceless property implies that the two Eigenvalues $\lambda_1,\ \lambda_2$ of elements of $\mathfrak{sl}(2,\mathbb{R})$ fulfil $\lambda_1+\lambda_2=0$. Let $A\in\mathfrak{sl}(2,\mathbb{R})\setminus\{0\}$ be arbitrary.
				
				Suppose first that $\lambda_1=\lambda_2=0$. By the assumption $A\ne 0$ we obtain that $A$ is under the adjoint action of $\mathrm{GL}(2,\mathbb{R})$ and, hence, $\mathrm{PGL}(2,\mathbb{R})$, equivalent to $N$. This follows from the fact that this is the only possible Jordan block for both Eigenvalues vanishing, such that $A$ itself is not zero.
				
				Next suppose that $\lambda_1$ is real and non-vanishing. Since $\lambda_2=-\lambda_1$, $A$ is equivalent to $|\lambda_1| E_3$, which immediately follows from the fact that $A$ has two Jordan blocks.
				
				Lastly, suppose that $\lambda_1\ne 0$ and that $\mathrm{Im}(\lambda_1)\ne 0$. Then $\lambda_2=\overline{\lambda_1}\ne\lambda_1$, and by the traceless property of $A$ we obtain $0=\lambda_1+\lambda_2=2\,\mathrm{Re}(\lambda_1)$. Hence, both $\lambda_1$ and $\lambda_2$ are purely imaginary, showing that $A$ is equivalent to $|\lambda_1|E_1$.
			\end{proof}
		\end{Lem}
	The action of $\mathrm{SL}(2,\mathbb{R})$ on the Lie algebra preserves the Killing form and, hence, the light cone and the level sets. In the proof of certain cases in Theorem \ref{thm_main} we will make use of the classification of symmetric $(0,2)$-tensors on $\mathbb{R}^3$ under $\mathrm{SO}(2,1)$. Note that in our convention, $\mathrm{SO}(2,1)$ leaves the bilinear form
		\begin{center}
			\scalebox{0.75}{
				$\left(\begin{array}{ccc}
					-1 & 0 & 0\\
					0 & 1 & 0\\
					0 & 0 & 1
				\end{array}
			\right)$}
		\end{center}
	invariant. Up to a change of basis, the following has been formulated as an exercise in \cite[pp.\,261,262]{O} and proven in \cite[Ch.\,5]{R}.
		
		\begin{Lem}\label{lem_sym2_mod_so21}
			Each class in the quotient space $\mathrm{Sym}^2(\mathbb{R}^3)^*/\mathrm{SO}(2,1)$ contains a representative of one of the five following forms:
				\begin{center}
					\begin{tabular}{ll}
						\scalebox{1}{
							$Q_1:=\left(\begin{array}{ccc}
									x & 0 & 0\\ 0 & y & 0\\ 0 & 0 & z
									\end{array}
								\right)$,$\vphantom{\begin{array}{c}0\\0\\0\\0\end{array}}$
							} &
						\scalebox{1}{
							$Q_2:=\left(\begin{array}{ccc}
									-x & y & 0\\ y & x & 0\\ 0 & 0 & z
									\end{array}
								\right)$,
							}\\
						\scalebox{1}{
							$Q_3:=\left(\begin{array}{ccc}
									-\frac{1}{2}-x & -\frac{1}{2} & 0\\ -\frac{1}{2} & -\frac{1}{2}+x & 0\\ 0 & 0 & y
									\end{array}
								\right)$,$\vphantom{\begin{array}{c}0\\0\\0\\0\end{array}}$
							} &
						\scalebox{1}{
							$Q_4:=\left(\begin{array}{ccc}
									\frac{1}{2}-x & \frac{1}{2} & 0\\ \frac{1}{2} & \frac{1}{2}+x & 0\\ 0 & 0 & y
									\end{array}
								\right)$,
							}\\
						\scalebox{1}{
							$Q_5:=\left(\begin{array}{ccc}
									-x & 0 & \frac{1}{\sqrt{2}}\\ 0 & x & \frac{1}{\sqrt{2}}\\ \frac{1}{\sqrt{2}} & \frac{1}{\sqrt{2}} & x
									\end{array}
								\right)$,$\vphantom{\begin{array}{c}0\\0\\0\\0\end{array}}$
							} &
					\end{tabular}
				\end{center}
				for some $x,y,z\in\mathbb{R}$.
		\end{Lem}
		
	Now we have all tools at hand to classify all one-parameter subgroups of $\mathrm{SL}(2,\mathbb{R})\times\mathrm{SL}(2,\mathbb{R})$ and the corresponding invariant scalar product on $\mathfrak{sl}(2,\mathbb{R})\times \mathfrak{sl}(2,\mathbb{R})$. To do so we will classify their infinitesimal generators up to equivalence under the induced action of the automorphism group.
		
		\begin{Prop}\label{prop_inf_gen_one_parameter_subgroups}
			The infinitesimal generator of a one-parameter subgroup of $\mathrm{SL}(2,\mathbb{R})\times\mathrm{SL}(2,\mathbb{R})$ is, up to equivalence, given by precisely one of the following:\newline
				\begin{minipage}[t]{0.4\textwidth}
					\begin{enumerate}[left=6em, label=(\roman*)]
						\item $(E_1,r E_1)$, $r\geq 0$,\label{eqn_generator_1}
						\item $(E_1,r E_3)$, $r> 0$,\label{eqn_generator_2}
						\item $(E_1,N)$,\label{eqn_generator_3}
					\end{enumerate}
				\end{minipage}
				\begin{minipage}[t]{0.4\textwidth}
					\begin{enumerate}[start=4, left=6em, label=(\roman*)]
						\item $(E_3,rE_3)$, $r\geq 0$,\label{eqn_generator_4}
						\item $(E_3,N)$,\label{eqn_generator_5}
						\item $(N,0)$,\label{eqn_generator_6}
						\item $(N,N)$.\label{eqn_generator_7}
					\end{enumerate}
				\end{minipage}
			\begin{proof}
				Two one-parameter subgroups coincide if their infinitesimal generators are related by a non-zero scaling factor. Hence, this proposition follows from Lemma \ref{lem_sl2_under_auto_class}, and by using the induced action of the $\mathbb{Z}_2$-part of the automorphism group of $\mathrm{SL}(2,\mathbb{R})\times\mathrm{SL}(2,\mathbb{R})$ which allows swapping the entries of $(A,B)\in\mathfrak{sl}(2,\mathbb{R})\oplus\mathfrak{sl}(2,\mathbb{R})$. Note that $N$ being lightlike with respect to the Killing form is the reason why we do not need to consider infinitesimal generators containing a scaling factor in front of the $N$-terms. This stems from the fact that the automorphism group of $\mathrm{SL}(2,\mathbb{R})$, that is $\mathrm{PGL}(2,\mathbb{R})$, acts with respect to the Killing form as $\mathrm{SO}(2,1)$ on $\mathfrak{sl}(2,\mathbb{R})$, and the action is transitive on the light cone.
			\end{proof}
		\end{Prop}
	
	Before employing the help of computer algebra software we still need to simplify the used input. To do so we classify all left-invariant pseudo-Riemannian metrics that are invariant under one of the one-parameter subgroups, or, equivalently, under their infinitesimal generators described in Proposition \ref{prop_inf_gen_one_parameter_subgroups}. Recall that a left-invariant pseudo-Riemannian metric is invariant under a one-parameter subgroup generated by $A\in\mathfrak{sl}(2,\mathbb{R})\oplus\mathfrak{sl}(2,\mathbb{R})$ if and only if $\mathrm{ad}_A$ is skew with respect to $g$, where $g$ as before denotes the value of the metric at the neutral element.
	
	\begin{Prop}\label{prop_classification_inv_subgrp_02tensors}
		Let $g$ be the value of a left-invariant symmetric $(0,2)$-tensor field on $\mathrm{SL}(2,\mathbb{R})\times\mathrm{SL}(2,\mathbb{R})$ at the neutral element. Then the tensor field is invariant under a one-parameter subgroup with infinitesimal generator one of Proposition \ref{prop_inf_gen_one_parameter_subgroups} \ref{eqn_generator_1}--\ref{eqn_generator_7} if and only if $g$ in the basis of $F_1,\ldots,F_6$ \eqref{eqn_basis_sl2_sl2} is of the form
			\begin{table}[H]\centering
					\begin{tabular}{|c|l|}
						\hline
						\text{$\overset{}{\vphantom{|}}$infinitesimal generator} & $\hspace*{0.21\textwidth}g$\\
						\hline\hline
						$(E_1,r E_1)$, $r>0$, $r\ne 1$ &\scalebox{0.75}{$\vphantom{\begin{array}{c} 0\\0\\0\\0\\0\\0\\0\end{array}}$
							$\left(\begin{array}{ccc|ccc}
										x_1 & 0 & 0 & a_1 & 0 & 0\\
										0 & y_1 & 0 & 0 & 0 & 0\\
										0 & 0 & y_1 & 0 & 0 & 0\\
										\hline
										\overset{}{\vphantom{|}}a_1 & 0 & 0 & x_2 & 0 & 0\\
										0 & 0 & 0 & 0 & y_2 & 0\\
										0 & 0 & 0 & 0 & 0 & y_2
									\end{array}
							\right)$}
						\\
						\hline
						$(E_1,E_1)$&\scalebox{0.75}{$\vphantom{\begin{array}{c} 0\\0\\0\\0\\0\\0\\0\end{array}}$
							$\left(\begin{array}{ccc|ccc}
										x_1 & 0 & 0 & a_1 & 0 & 0\\
										0 & y_1 & 0 & 0 & b_2 & c_2\\
										0 & 0 & y_1 & 0 & -c_2 & b_2\\
										\hline
										\overset{}{\vphantom{|}}a_1 & 0 & 0 & x_2 & 0 & 0\\
										0 & b_2 & -c_2 & 0 & y_2 & 0\\
										0 & c_2 & b_2 & 0 & 0 & y_2
									\end{array}
							\right)$}
						\\
						\hline
						$(E_1,0)$&\scalebox{0.75}{$\vphantom{\begin{array}{c} 0\\0\\0\\0\\0\\0\\0\end{array}}$
							$\left(\begin{array}{ccc|ccc}
										x_1 & 0 & 0 & a_1 & b_1 & c_1\\
										0 & y_1 & 0 & 0 & 0 & 0\\
										0 & 0 & y_1 & 0 & 0 & 0\\
										\hline
										\overset{}{\vphantom{|}}a_1 & 0 & 0 &  &  & \\
										b_1 & 0 & 0 &  & Q_i & \\
										c_1 & 0 & 0 &  &  & 
									\end{array}
							\right)$, $1\leq i\leq 5$}
						\\
						\hline
						$(E_1,r E_3)$, $r> 0$, $r\ne 1$ &\scalebox{0.75}{$\vphantom{\begin{array}{c} 0\\0\\0\\0\\0\\0\\0\end{array}}$
							$\left(\begin{array}{ccc|ccc}
										x_1 & 0 & 0 & 0 & 0 & c_1\\
										0 & y_1 & 0 & 0 & 0 & 0\\
										0 & 0 & y_1 & 0 & 0 & 0\\
										\hline
										\overset{}{\vphantom{|}}0 & 0 & 0 & x_2 & 0 & 0\\
										0 & 0 & 0 & 0 & -x_2 & 0\\
										c_1 & 0 & 0 & 0 & 0 & z_2
									\end{array}
							\right)$}
						\\
						\hline
						$(E_1,N)$ &\scalebox{0.75}{$\vphantom{\begin{array}{c} 0\\0\\0\\0\\0\\0\\0\end{array}}$
							$\left(\begin{array}{ccc|ccc}
										x_1 & 0 & 0 & a_1 & a_1 & 0\\
										0 & y_1 & 0 & 0 & 0 & 0\\
										0 & 0 & y_1 & 0 & 0 & 0\\
										\hline
										\overset{}{\vphantom{|}}a_1 & 0 & 0 & x_2 & (x_2+y_2)/2 & 0\\
										a_1 & 0 & 0 & (x_2+y_2)/2 & y_2 & 0\\
										0 & 0 & 0 & 0 & 0 & (-x_2+y_2)/2
									\end{array}
							\right)$}
						\\
						\hline
						$(E_3,rE_3)$, $r>0$, $r\ne 1$ &\scalebox{0.75}{$\vphantom{\begin{array}{c} 0\\0\\0\\0\\0\\0\\0\end{array}}$
							$\left(\begin{array}{ccc|ccc}
										x_1 & 0 & 0 & 0 & 0 & 0\\
										0 & -x_1 & 0 & 0 & 0 & 0\\
										0 & 0 & z_1 & 0 & 0 & c_3\\
										\hline
										\overset{}{\vphantom{|}}0 & 0 & 0 & x_2 & 0 & 0\\
										0 & 0 & 0 & 0 & -x_2 & 0\\
										0 & 0 & c_3 & 0 & 0 & z_2
									\end{array}
							\right)$}
						\\
						\hline
					\end{tabular}
			\caption{}\label{table_infgen_g_1}
			\end{table}
			\begin{table}[H]\centering
					\begin{tabular}{|c|l|}
						\hline
						\text{$\overset{}{\vphantom{|}}$infinitesimal generator} & $\hspace*{0.33\textwidth}g$\\
						\hline\hline
						$(E_3,E_3)$&\scalebox{0.75}{$\vphantom{\begin{array}{c} 0\\0\\0\\0\\0\\0\\0\end{array}}$
							$\left(\begin{array}{ccc|ccc}
										x_1 & 0 & 0 & a_1 & b_1 & 0\\
										0 & -x_1 & 0 & -b_1 & -a_1 & 0\\
										0 & 0 & z_1 & 0 & 0 & c_3\\
										\hline
										\overset{}{\vphantom{|}}a_1 & -b_1 & 0 & x_2 & 0 & 0\\
										b_1 & -a_1 & 0 & 0 & -x_2 & 0\\
										0 & 0 & c_3 & 0 & 0 & z_2
									\end{array}
							\right)$}
						\\
						\hline
						$(E_3,0)$&\scalebox{0.75}{$\vphantom{\begin{array}{c} 0\\0\\0\\0\\0\\0\\0\end{array}}$
							$\left(\begin{array}{ccc|ccc}
										x_1 & 0 & 0 & 0 & 0 & 0\\
										0 & -x_1 & 0 & 0 & 0 & 0\\
										0 & 0 & z_1 & a_3 & b_3 & c_3\\
										\hline
										\overset{}{\vphantom{|}}0 & 0 & a_3 &  &  & \\
										0 & 0 & b_3 &  & Q_i & \\
										0 & 0 & c_3 &  &  & 
									\end{array}
							\right)$, $1\leq i\leq 5$}
						\\
						\hline
						$(E_3,N)$ &\scalebox{0.75}{$\vphantom{\begin{array}{c} 0\\0\\0\\0\\0\\0\\0\end{array}}$
							$\left(\begin{array}{ccc|ccc}
										x_1 & 0 & 0 & 0 & 0 & 0\\
										0 & -x_1 & 0 & 0 & 0 & 0\\
										0 & 0 & z_1 & a_3 & a_3 & 0\\
										\hline
										\overset{}{\vphantom{|}}0 & 0 & a_3 & x_2 & (x_2+y_2)/2 & 0\\
										0 & 0 & a_3 & (x_2+y_2)/2 & y_2 & 0\\
										0 & 0 & 0 & 0 & 0 & (-x_2+y_2)/2
									\end{array}
							\right)$}
						\\
						\hline
						$(N,0)$ &\scalebox{0.75}{$\vphantom{\begin{array}{c} 0\\0\\0\\0\\0\\0\\0\end{array}}$
							$\left(\begin{array}{ccc|ccc}
										x_1 & (x_1+y_1)/2 & 0 & a_1 & b_1 & c_1\\
										(x_1+y_1)/2 & y_1 & 0 & a_1 & b_1 & c_1\\
										0 & 0 & (-x_1+y_1)/2 & 0 & 0 & 0\\
										\hline
										\overset{}{\vphantom{|}}a_1 & a_1 & 0 &  &  & \\
										b_1 & b_1 & 0 &  & Q_i & \\
										c_1 & c_1 & 0 &  &  & 
									\end{array}
							\right)$, $1\leq i\leq 5$}
						\\
						\hline
						$(N,N)$ &\scalebox{0.75}{$\vphantom{\begin{array}{c} 0\\0\\0\\0\\0\\0\\0\end{array}}$
							$\left(\begin{array}{ccc|ccc}
										x_1 & (x_1+y_1)/2 & 0 & a_1 & b_1 & c_1\\
										(x_1+y_1)/2 & y_1 & 0 & b_1 & -a_1+2b_1 & c_1\\
										0 & 0 & (-x_1+y_1)/2 & -c_1 & -c_1 & -a_1+b_1\\
										\hline
										\overset{}{\vphantom{|}}a_1 & b_1 & -c_1 & x_2 & (x_2+y_2)/2 & 0 \\
										b_1 & -a_1+2b_1 & -c_1 & (x_2+y_2)/2 & y_2 & 0 \\
										c_1 & c_1 & -a_1+b_1 & 0 & 0 & (-x_2+y_2)/2
									\end{array}
							\right)$}
						\\
						\hline
					\end{tabular}					
			\caption{}\label{table_infgen_g_2}
			\end{table}
		\noindent
		for some $x_1,y_1,z_1,x_2,y_2,z_2,a_1,b_1,c_1,a_2,b_2,c_2,a_3,b_3,c_3\in\mathbb{R}$ with $Q_i$, $1\leq i\leq 5$, as in Lemma \ref{lem_sym2_mod_so21}.
		\begin{proof}
			First we will classify all left-invariant symmetric $(0,2)$-tensors on $\mathrm{SL}(2,\mathbb{R})$ that are invariant under $E_1$, $E_3$, or $N$, respectively. In the basis $E_1,E_2,E_3$, a left-invariant symmetric $(0,2)$-tensors on $\mathrm{SL}(2,\mathbb{R})$ at the neutral element is of the form
				\begin{equation*}
					g=\left(\begin{array}{ccc}
								x & u & w\\
								u & y & v\\
								w & v & z
							\end{array}
					\right),\quad x,y,z,u,v,w\in\mathbb{R},
				\end{equation*}
			and it is invariant under $A\in\mathfrak{sl}(2,\mathbb{R})$ if and only if $\mathrm{ad}_A^T\,g + g\,\mathrm{ad}_A=0$. We find that $g$ is invariant under $A\in\{E_1,E_3,N\}$ if and only if
				\begin{table}[H]\centering
				\begin{tabular}{|c|l|}
					\hline
					$\overset{}{\vphantom{|}} A$ &\hspace*{0.14\textwidth} $g$\\
					\hline\hline
					$E_1$ & \scalebox{0.75}{
							$\vphantom{\begin{array}{c} 0\\0\\0\\0\end{array}}\left(\begin{array}{ccc}
										x & 0 & 0\\
										0 & y & 0\\
										0 & 0 & y
								\end{array}
							\right)$}\\
					\hline
					$E_3$ & \scalebox{0.75}{
							$\vphantom{\begin{array}{c} 0\\0\\0\\0\end{array}}\left(\begin{array}{ccc}
										x & 0 & 0\\
										0 & -x & 0\\
										0 & 0 & z
								\end{array}
							\right)$}\\
					\hline
					$N$ & \scalebox{0.75}{
							$\vphantom{\begin{array}{c} 0\\0\\0\\0\end{array}}\left(\begin{array}{ccc}
										x & (x+y)/2 & 0\\
										(x+y)/2 & y & 0\\
										0 & 0 & (-x+y)/2
								\end{array}
							\right)$}\\
					\hline
				\end{tabular}
				\end{table}
			\noindent
			for some $x,y,z\in\mathbb{R}$. Together with Lemma \ref{lem_sym2_mod_so21} this already shows that the diagonal blocks in Tables \ref{table_infgen_g_1} and \ref{table_infgen_g_2} are correct. For the off-diagonal blocks we first note that a symmetric $(0,2)$-tensor on $\mathrm{SL}(2,\mathbb{R})\times\mathrm{SL}(2,\mathbb{R})$ at the neutral element in the basis $F_1,\ldots,F_6$ is of the form
				\begin{equation}
					g=\left(\begin{array}{ccc|ccc}
								x_1 & u_1 & w_1 & a_1 & b_1 & c_1\\
								u_1 & y_1 & v_1 & a_2 & b_2 & c_2\\
								w_1 & v_1 & z_1 & a_3 & b_3 & c_3\\
								\hline
								a_1 & a_2 & a_3 & x_2 & u_2 & w_2\\
								b_1 & b_2 & b_3 & u_2 & y_2 & v_2\\
								c_1 & c_2 & c_3 & w_2 & v_2 & z_2
							\end{array}
					\right)\label{eqn_g_sl2sl2_general_form}
				\end{equation}
			with real entries. For the off-diagonal blocks to be invariant under $(A_1,A_2)\in\mathfrak{sl}(2,\mathbb{R})\oplus\mathfrak{sl}(2,\mathbb{R})$, they need to fulfil
				\begin{equation*}
					\mathrm{ad}_B^T\,\left(\begin{array}{ccc}a_1 & a_2 & a_3\\b_1 & b_2 & b_3\\c_1 & c_2 & c_3\end{array}\right)+\left(\begin{array}{ccc}a_1 & a_2 & a_3\\b_1 & b_2 & b_3\\c_1 & c_2 & c_3\end{array}\right)\,\mathrm{ad}_A=0.
				\end{equation*}
			Now a straight-forward calculation with $A$, $B$, matching the cases of the infinitesimal generators in Tables \ref{table_infgen_g_1} and \ref{table_infgen_g_2} leads to the claimed results.
		\end{proof}
	\end{Prop}
	
\section{Classification result}
	The next step is to implement the results of Proposition \ref{prop_classification_inv_subgrp_02tensors} in a computer algebra system to solve the Einstein equation $\mathrm{Ric}=\lambda\,g$ for a left-invariant pseudo-Riemannian metric $g$ on $\mathrm{SL}(2,\mathbb{R})\times\mathrm{SL}(2,\mathbb{R})$ that is additionally invariant under a one-parameter subgroup. To do so we use the formulas for the connection coefficients \eqref{eqn_con_coeff} respectively the components of the Riemann curvature tensor \eqref{eqn_riemann_curvature} discussed in the introduction. As before, $g$ will now denote the value of the considered metric at the neutral element.

	For the computer algebra system we use Maple. To execute our script to find all such pseudo-Riemannian Einstein metrics, at this point not necessarily up to isometry, simply run \href{http://www.math.uni-hamburg.de/home/cortes/sl2_sl2_einstein.txt}{\texttt{maple sl2\char`_sl2\char`_einstein.txt}} in any Unix-like environment with an installation of Maple. See the provided commentary for different output options. It should run correctly on Maple 2019 or later. The general idea is to calculate an explicit expression of the equation $\mathrm{Ric}-\lambda\,g=0$ at the neutral element and obtain $21$ equations that are, in general, fractions of polynomials in the variables $x_1,\ldots,c_3$ as in \eqref{eqn_g_sl2sl2_general_form}. Since rescaling a pseudo-Einstein metric with Einstein constant $\lambda$ with a non-zero constant $C$ yields another pseudo-Einstein metric with Einstein constant $C^2$ we will also require $\det(g)=\pm 1$. Here $\det(g)$ is to be understood as the determinant of the right-hand side of \eqref{eqn_g_sl2sl2_general_form}, which can be viewed as the volume of the pseudo-Riemannian metric at the neutral element with respect to the chosen basis $F_1,\ldots,F_6$. All in all, for each sign of $\det(g)$ and each case listed in Tables \ref{table_infgen_g_1} and \ref{table_infgen_g_2} we have $22$ equations that are all fractions of polynomials and that we all need to solve in order to obtain the presumed pseudo-Riemannian Einstein metrics.

	We now take the numerators of said equations and calculate in each case a lexicographic Gr\"obner basis in the variable order
		\begin{equation*}
			\lambda,x_1,y_1,z_1,x_2,y_2,z_2,u_1,v_1,w_1,u_2,v_2,w_2,a_1,b_1,c_1,a_2,b_2,c_2,a_3,b_3,c_3.
		\end{equation*}
	This variable order is not in any way canonical, but after some testing it worked decently fast in every single case, see the comments in the our code for details of expected RAM use and approximate wall time. We then set the so-obtained entries of the Gr\"obner basis to zero and obtain potential solutions of the Einstein equation in each case, ``potential'' since they are only solutions of the numerator of $\mathrm{Ric}-\lambda\,g=0$ at the neutral element.

	Lastly we take the so-obtained potential solutions and, by evaluating $\mathrm{Ric}-\lambda\,g$ at the neutral element, verify that in each case they are actual solutions, thus yielding all possible left-invariant pseudo-Riemannian Einstein metrics on $\mathrm{SL}(2,\mathbb{R})\times\mathrm{SL}(2,\mathbb{R})$ that are invariant under a one-parameter subgroup. In the step thereafter, we will take the results and check by hand which of the solutions are isometric.
		
		\begin{rem}
			If $g$ is a solution of the Einstein equation with Einstein constant $\lambda$ on an even-dimensional manifold, then $-g$ is also a pseudo-Riemannian Einstein metric with Einstein constant $-\lambda$.
		\end{rem}
		
		The above remark and $\dim(\mathrm{SL}(2,\mathbb{R})\times\mathrm{SL}(2,\mathbb{R}))=6$ imply that in our calculations we can without loss of generality assume that the Einstein constant is non-positive. 
		
		\begin{Lem}\label{lem_isometry_class_solutions}
			Let $g$ denote the value of a solution of the Einstein equation at the neutral element of $\mathrm{SL}(2,\mathbb{R})\times\mathrm{SL}(2,\mathbb{R})$ with non-positive Einstein constant $\lambda$ that is invariant under a one-parameter subgroup with infinitesimal generator one of Proposition \ref{prop_inf_gen_one_parameter_subgroups} \ref{eqn_generator_1}--\ref{eqn_generator_6}. Let $B$ denote the Killing form of $\mathrm{SL}(2,\mathbb{R})\times\mathrm{SL}(2,\mathbb{R})$. Then $g$ is given by one of the following entries of Tables \ref{table_solutions_positive_lambda_not_mod_isom_1} and \ref{table_solutions_positive_lambda_not_mod_isom_2}:
			\begin{table}[H]\centering
					\begin{tabular}{|c|l|}
						\hline
						\text{$\overset{}{\vphantom{|}}$infinitesimal generator} & solutions of $\mathrm{Ric}-\lambda\,g=0$ with $\lambda\leq 0$ at the neutral element\\
						\hline\hline
						$(E_1,r E_1)$, $r>0$, $r\ne 1$ &\scalebox{0.75}{$\vphantom{\begin{array}{c} 0\\0\\0\\0\\0\\0\\0\end{array}}$
							$g=\frac{1}{8}B=\left(\begin{array}{ccc|ccc}
										-1 & 0 & 0 & 0 & 0 & 0\\
										0 & 1 & 0 & 0 & 0 & 0\\
										0 & 0 & 1 & 0 & 0 & 0\\
										\hline
										\overset{}{\vphantom{|}}0 & 0 & 0 & -1 & 0 & 0\\
										0 & 0 & 0 & 0 & 1 & 0\\
										0 & 0 & 0 & 0 & 0 & 1
									\end{array}
							\right)$, $\lambda=-2$, $\det(g)=1$, $\mathrm{sgn}(g)=2$}
						\\
						\hline
						\multirow{10.5}{*}
						{$(E_1,E_1)$}&\scalebox{0.75}{$\vphantom{\begin{array}{c} 0\\0\\0\\0\\0\\0\\0\\0\end{array}}$
							$g=\left(\begin{array}{ccc|ccc}
										-(1+k) & 0 & 0 & 1 & 0 & 0\\
										0 & 1+k & 0 & 0 & \pm\sqrt{1-c_2^2} & c_2\\
										0 & 0 & 1+k & 0 & -c_2 & \pm\sqrt{1-c_2^2}\\
										\hline
										\overset{}{\vphantom{|}}1 & 0 & 0 & -(2-k) & 0 & 0\\
										0 & \pm\sqrt{1-c_2^2} & -c_2 & 0 & 2-k & 0\\
										0 & c_2 & \pm\sqrt{1-c_2^2} & 0 & 0 & 2-k
									\end{array}
							\right)$,}
							\\
							&\scalebox{0.75}{$k\in\{0,1\}$, $|c_2|\leq 1$, $\lambda=-2$, $\det(g)=1$, $\mathrm{sgn}(g)=2$ $\vphantom{\begin{array}{c} 0\\0\end{array}}$}
						\\ \cline{2-2}
							&\scalebox{0.75}{$\vphantom{\begin{array}{c} 0\\0\\0\\0\\0\\0\\0\\0\end{array}}$
								$g=\left(\begin{array}{ccc|ccc}
											-2/\sqrt{3} & 0 & 0 & 1/\sqrt{3} & 0 & 0\\
											0 & 2/\sqrt{3} & 0 & 0 & \pm\sqrt{1/3-c_2^2} & c_2\\
											0 & 0 & 2/\sqrt{3} & 0 & -c_2 & \pm\sqrt{1/3-c_2^2}\\
											\hline
											\overset{}{\vphantom{|}}1/\sqrt{3} & 0 & 0 & -2/\sqrt{3} & 0 & 0\\
											0 & \pm\sqrt{1/3-c_2^2} & -c_2 & 0 & 2/\sqrt{3} & 0\\
											0 & c_2 & \pm\sqrt{1/3-c_2^2} & 0 & 0 & 2/\sqrt{3}
										\end{array}
								\right)$,}
								\\
								&\scalebox{0.75}{$k\in\{0,1\}$, $|c_2|\leq 1/\sqrt{3}$, $\lambda=-\frac{10}{3\sqrt{3}}$, $\det(g)=1$, $\mathrm{sgn}(g)=2$ $\vphantom{\begin{array}{c} 0\\0\end{array}}$}
						\\ \cline{2-2}
							&\scalebox{0.75}{$\vphantom{\begin{array}{c} 0\\0\end{array}}$
								$g=\frac{1}{8}B$, $\lambda=-2$, $\det(g)=1$, $\mathrm{sgn}(g)=2$}
						\\
						\hline
						$(E_1,0)$&\scalebox{0.75}{$\vphantom{\begin{array}{c} 0\\0\end{array}}$
								$g=\frac{1}{8}B$, $\lambda=-2$, $\det(g)=1$, $\mathrm{sgn}(g)=2$}
						\\
						\hline
						$(E_1,r E_3)$, $r> 0$, $r\ne 1$ &\scalebox{0.75}{$\vphantom{\begin{array}{c}0\\0\end{array}}$
								$g=\frac{1}{8}B$, $\lambda=-2$, $\det(g)=1$, $\mathrm{sgn}(g)=2$}
						\\
						\hline
						$(E_1,N)$ &\scalebox{0.75}{$\vphantom{\begin{array}{c} 0\\0\end{array}}$
								$g=\frac{1}{8}B$, $\lambda=-2$, $\det(g)=1$, $\mathrm{sgn}(g)=2$}
						\\
						\hline
						$(E_3,rE_3)$, $r>0$, $r\ne 1$ &\scalebox{0.75}{$\vphantom{\begin{array}{c}0\\0\end{array}}$
								$g=\frac{1}{8}B$, $\lambda=-2$, $\det(g)=1$, $\mathrm{sgn}(g)=2$}
						\\
						\hline
						\multirow{10.5}{*}{$(E_3,E_3)$}&\scalebox{0.75}{$\vphantom{\begin{array}{c} 0\\0\\0\\0\\0\\0\\0\\0\end{array}}$
							$g=\left(\begin{array}{ccc|ccc}
										-(1+k) & 0 & 0 & \pm\sqrt{1+b_1^2} & b_1 & 0\\
										0 & (1+k) & 0 & -b_1 & \mp\sqrt{1+b_1^2} & 0\\
										0 & 0 & (1+k) & 0 & 0 & -1\\
										\hline
										\overset{}{\vphantom{\big|^|}}\pm\sqrt{1+b_1^2} & -b_1 & 0 & -(2-k) & 0 & 0\\
										b_1 & \mp\sqrt{1+b_1^2} & 0 & 0 & 2-k & 0\\
										0 & 0 & -1 & 0 & 0 & 2-k
									\end{array}
							\right)$,}
						\\
							&\scalebox{0.75}{$k\in\{0,1\}$, $b_1\in\mathbb{R}$, $\lambda=-2$, $\det(g)=1$, $\mathrm{sgn}(g)=2$ $\vphantom{\begin{array}{c} 0\\0\end{array}}$}
						\\ \cline{2-2}
							&\scalebox{0.75}{$\vphantom{\begin{array}{c} 0\\0\\0\\0\\0\\0\\0\\0\end{array}}$
							$g=\left(\begin{array}{ccc|ccc}
										-2/\sqrt{3} & 0 & 0 & \pm\sqrt{1/3+b_1^2} & b_1 & 0\\
										0 & 2/\sqrt{3} & 0 & -b_1 & \mp\sqrt{1/3+b_1^2} & 0\\
										0 & 0 & 2/\sqrt{3} & 0 & 0 & -1/\sqrt{3} \\
										\hline
										\overset{}{\vphantom{\big|^|}}\pm\sqrt{1/3+b_1^2} & -b_1 & 0 & -2/\sqrt{3} & 0 & 0\\
										b_1 & \mp\sqrt{1/3+b_1^2} & 0 & 0 & 2/\sqrt{3} & 0\\
										0 & 0 & -1/\sqrt{3} & 0 & 0 & 2/\sqrt{3}
									\end{array}
							\right)$,}
						\\
							&\scalebox{0.75}{$b_1\in\mathbb{R}$, $\lambda=-\frac{10}{3\sqrt{3}}$, $\det(g)=1$, $\mathrm{sgn}(g)=2$ $\vphantom{\begin{array}{c} 0\\0\end{array}}$}
						\\ \cline{2-2}
							&\scalebox{0.75}{$\vphantom{\begin{array}{c}0\\0\end{array}}$
								$g=\frac{1}{8}B$, $\lambda=-2$, $\det(g)=1$, $\mathrm{sgn}(g)=2$}
						\\
						\hline
					\end{tabular}
			\caption{}\label{table_solutions_positive_lambda_not_mod_isom_1}
			\end{table}
			\begin{table}[H]\centering
				\begin{tabular}{|c|l|}
					\hline
					\text{$\overset{}{\vphantom{|}}$infinitesimal generator} & solutions of $\mathrm{Ric}-\lambda\,g=0$ with $\lambda\leq 0$ at the neutral element\\
					\hline\hline
					$(E_3,0)$&\scalebox{0.75}{$\vphantom{\begin{array}{c} 0\\0\end{array}}$
								$g=\frac{1}{8}B$, $\lambda=-2$, $\det(g)=1$, $\mathrm{sgn}(g)=2$}
					\\
					\hline
					$(E_3,N)$ &\scalebox{0.75}{$\vphantom{\begin{array}{c} 0\\0\end{array}}$
									$g=\frac{1}{8}B$, $\lambda=-2$, $\det(g)=1$, $\mathrm{sgn}(g)=2$}
					\\
					\hline
					$(N,0)$ &\scalebox{0.75}{$\vphantom{\begin{array}{c} 0\\0\end{array}}$
									$g=\frac{1}{8}B$, $\lambda=-2$, $\det(g)=1$, $\mathrm{sgn}(g)=2$}
					\\
					\hline
					\multirow{10.5}{*}{$(N,N)$} &\scalebox{0.75}{$\vphantom{\begin{array}{c} 0\\0\\0\\0\\0\\0\\0\end{array}}$
						$g=\left(\begin{array}{ccc|ccc}
									-(1+k) & 0 & 0 & 1+c_1^2/2 & c_1^2/2 & c_1\\
									0 & 1+k & 0 & c_1^2/2 & c_1^2/2-1 & c_1\\
									0 & 0 & 1+k & -c_1 & -c_1 & -1\\
									\hline
									\overset{}{\vphantom{|}}1+c_1^2/2 & c_1^2/2 & -c_1 & -(2-k) & 0 & 0 \\
									c_1^2/2 & c_1^2/2-1 & -c_1 & 0 & 2-k & 0 \\
									c_1 & c_1 & -1 & 0 & 0 & 2-k
								\end{array}
						\right)$,}
					\\
						&\scalebox{0.75}{$k\in\{0,1\}$, $c_1\in\mathbb{R}$, $\lambda=-2$, $\det(g)=1$, $\mathrm{sgn}(g)=2$ $\vphantom{\begin{array}{c} 0\\0\end{array}}$}
					\\ \cline{2-2}
						&\scalebox{0.75}{$\vphantom{\begin{array}{c} 0\\0\\0\\0\\0\\0\\0\end{array}}$
						$g=\left(\begin{array}{ccc|ccc}
									-2/\sqrt{3} & 0 & 0 & \sqrt{3}(2+3c_1^2)/6 & \sqrt{3}c_1^2/2 & c_1\\
									0 & 2/\sqrt{3} & 0 & \sqrt{3}c_1^2/2 & \sqrt{3}(-2+3c_1^2)/6 & c_1\\
									0 & 0 & 2/\sqrt{3} & -c_1 & -c_1 & -1/\sqrt{3} \\
									\hline
									\overset{}{\vphantom{|}}\sqrt{3}(2+3c_1^2)/6 & \sqrt{3}c_1^2/2 & -c_1 & -2/\sqrt{3} & 0 & 0\\
									\sqrt{3}c_1^2/2 & \sqrt{3}(-2+3c_1^2)/6 & -c_1 & 0 & 2/\sqrt{3} & 0\\
									c_1 & c_1 & -1/\sqrt{3} & 0 & 0 & 2/\sqrt{3}
								\end{array}
						\right)$,}
					\\
						&\scalebox{0.75}{$c_1\in\mathbb{R}$, $\lambda=-\frac{10}{3\sqrt{3}}$, $\det(g)=1$, $\mathrm{sgn}(g)=2$ $\vphantom{\begin{array}{c} 0\\0\end{array}}$}
					\\ \cline{2-2}
						&\scalebox{0.75}{$\vphantom{\begin{array}{c} 0\\0\end{array}}$
								$g=\frac{1}{8}B$, $\lambda=-2$, $\det(g)=1$, $\mathrm{sgn}(g)=2$}
					\\
					\hline
				\end{tabular}					
			\caption{}\label{table_solutions_positive_lambda_not_mod_isom_2}
			\end{table}
			\begin{proof}
				These are precisely the results obtained via Maple.
			\end{proof}
		\end{Lem}
		
		In the next step of the proof of Theorem \ref{thm_main} we need to determine the equivalence classes of the above solutions modulo the automorphisms of $\mathrm{SL}(2,\mathbb{R})\times\mathrm{SL}(2,\mathbb{R})$.
		
		\begin{Lem}\label{lem_family_solutions_representative}
			Under the automorphism group $(\mathrm{PGL}(2,\mathbb{R})\times\mathrm{PGL}(2,\mathbb{R}))\rtimes\mathbb{Z}_2$, each family of solutions in Lemma \ref{lem_isometry_class_solutions} belongs to a single orbit.
			\begin{proof}
				We define
					\begin{center}
						\begin{tabular}{ll}
							\scalebox{0.75}{$g_1:=
								\left(\begin{array}{ccc|ccc}
											-1 & 0 & 0 & 1 & 0 & 0\\
											0 & 1 & 0 & 0 & 1 & 0\\
											0 & 0 & 1 & 0 & 0 & 1\\
											\hline
											\overset{}{\vphantom{|}}1 & 0 & 0 & -2 & 0 & 0\\
											0 & 1 & 0 & 0 & 2 & 0\\
											0 & 0 & 1 & 0 & 0 & 2
										\end{array}
								\right)$,}
							&
							\scalebox{0.75}{$g_2:=
								\frac{1}{\sqrt{3}}\left(\begin{array}{ccc|ccc}
											-2 & 0 & 0 & 1 & 0 & 0\\
											0 & 2 & 0 & 0 & 1 & 0\\
											0 & 0 & 2 & 0 & 0 & 1\\
											\hline
											\overset{}{\vphantom{|}}1 & 0 & 0 & -2 & 0 & 0\\
											0 & 1 & 0 & 0 & 2 & 0\\
											0 & 0 & 1 & 0 & 0 & 2
										\end{array}
								\right)$.}
						\end{tabular}
					\end{center}
				For the first family of solutions corresponding to $(E_1,E_1)$ in Table \ref{table_solutions_positive_lambda_not_mod_isom_1} we can set $k=0$ after a possible swap of the $\mathrm{SL}(2,\mathbb{R})$-factors. We see that
					\begin{center}
						\scalebox{0.75}{$
						\left(
							\begin{array}{cc}
								\pm\sqrt{1/3-c_2^2} & -c_2 \\
								c_2 & \pm\sqrt{1/3-c_2^2} 
							\end{array}
						\right)$}
					\end{center}
				is a special orthogonal matrix for all allowed values of $c_1$. By the invariance of the upper left block or, equivalently, lower right block of the solutions under $\{1\}\times\mathrm{SO(2)}\subset\mathrm{SO}(2,1)$ it follows after a possible orthogonal transformation in $x_1,y_1$, respectively $x_2,y_2$, that each element of the family of solutions is, under the outer automorphism group of $\mathrm{SL}(2,\mathbb{R})\times\mathrm{SL}(2,\mathbb{R})$, equivalent to $g_1$. The argument for the second family of solutions corresponding to $(E_1,E_1)$ is analogous and we obtain that each element is equivalent to $g_2$.

				The case of the two families of solutions corresponding to $(E_3,E_3)$ in Table \ref{table_solutions_positive_lambda_not_mod_isom_1} works similarly with the only difference that
					\begin{center}
						\scalebox{0.75}{$
						\left(
							\begin{array}{cc}
								\pm\sqrt{1+b_1^2} & -b_1\\
								b_1 & \mp\sqrt{1+b_1^2}
							\end{array}
						\right)$}
					\end{center}
				and the corresponding matrix in the second family of solutions is now not special orthogonal, but rather an element in $\mathrm{O}(1,1)$. Note that while the determinant of the above matrix is negative, the entire matrix in the lower left block is still contained in $\mathrm{SO}(2,1)$ because of the only other non-zero entry $-1/\sqrt{3}$. This also holds for the second family of solutions. For the first family of solutions we find that every element is equivalent to $g_1$ and the elements of the second family of solutions are equivalent to $g_2$.

				Lastly, we have to deal with the families of solutions corresponding to the infinitesimal generator $(N,N)$. In this case it is a little more difficult to find the correct transformations in the automorphism group, but the argument works in the same way as for the other families of solutions. Let
					\begin{center}
						\begin{tabular}{ll}
							\scalebox{0.75}{$T_1:=
								\left(\begin{array}{ccc}
											1+c_1^2/2 & c_1^2/2 & -c_1 \\
											-c_1^2/2 & 1-c_1^2/2 & c_1\\
											-c_1 & -c_1 & 1 
										\end{array}
								\right)$,}
							&
							\scalebox{0.75}{$T_2:=
								\left(\begin{array}{ccc}
											1+3c_1^2/2 & 3c_1^2/2 & -\sqrt{3}c_1 \\
											-3c_1^2/2 & 1-3c_1^2/2 & \sqrt{3}c_1\\
											-\sqrt{3}c_1 & -\sqrt{3}c_1 & 1 
										\end{array}
								\right)$,}
						\end{tabular}
					\end{center}
				with $c_1\in\mathbb{R}$. We check that both $T_1$ and $T_2$ are contained in $\mathrm{SO}(2,1)$ for all $c_1\in\mathbb{R}$. For the first family of solutions corresponding to $(N,N)$ we transform $x_1,y_1,z_1$ with $T_1$ and calculate
					\begin{center}
						\scalebox{0.75}{$T_1\cdot
							\left(\begin{array}{ccc}
										1+c_1^2/2 & c_1^2/2 & -c_1 \\
										c_1^2/2 & c_1^2/2-1 & -c_1 \\
										c_1 & c_1 & -1 
									\end{array}
							\right)=\left(\begin{array}{ccc} 1& 0 & 0\\ 0 & 1 & 0\\ 0 & 0 & 1\end{array}\right)$,}
					\end{center}
				and for the second family of solutions we proceed analogously and obtain
					\begin{center}
						\scalebox{0.75}{$T_2\cdot
							\left(\begin{array}{ccc}
										\sqrt{3}(2+3c_1^2)/6 & \sqrt{3}c_1^2/2 & -c_1 \\
										\sqrt{3}c_1^2/2 & \sqrt{3}(-2+3c_1^2)/6 & -c_1 \\
										c_1 & c_1 & -1/\sqrt{3} 
									\end{array}
							\right)=\frac{1}{\sqrt{3}}\left(\begin{array}{ccc} 1& 0 & 0\\ 0 & 1 & 0\\ 0 & 0 & 1\end{array}\right)$.}
					\end{center}
				Hence, every element of the first family of solutions is equivalent to $g_1$, and every element of the second family of solutions is equivalent to $g_2$.
			\end{proof}
		\end{Lem}

		In order to complete the proof of Theorem \ref{thm_main}, we need to verify that $\frac{1}{8}B$ is isometric to $g_1$, but not isometric to $g_2$.

		\begin{Lem}\label{lem_equiv_B_g1}
			The left-invariant metrics corresponding to $\frac{1}{8}B$ and $g_1$ are isometric.
			\begin{proof}
				The validity of this lemma is not immediately clear, as it turns out that $\frac{1}{8}B$ and $g_1$ are not isometric via a transformation in $\mathrm{Aut}(\mathrm{SL}(2,\mathbb{R})\times\mathrm{SL}(2,\mathbb{R}))$. Using our Maple script \href{http://www.math.uni-hamburg.de/home/cortes/ssli_cov_R.txt}{\texttt{ssli\char`_cov\char`_R.txt}}, respectively its application in \href{http://www.math.uni-hamburg.de/home/cortes/ssli_cov_R_Killing_form.txt}{\texttt{ssli\char`_cov\char`_R\char`_Killing\char`_form.txt}} and \href{http://www.math.uni-hamburg.de/home/cortes/ssli_cov_R_g1.txt}{\texttt{ssli\char`_cov\char`_R\char`_g1.txt}}, we find that both $\nabla^{\frac{1}{8}B}R^{\frac{1}{8}B}$ and $\nabla^{g_1}R^{g_1}$ vanish identically. Note that the first claim automatically follows from the fact that $\mathrm{SL}(2,\mathbb{R})\times\mathrm{SL}(2,\mathbb{R})$ equipped with its Killing form is a symmetric space. By a generalisation of a result from \cite{AW}, cf. \cite[Sect.\,5]{W}, showing that $\frac{1}{8}B$ and $g_1$ are isometric amounts to finding a linear transformation $A:\mathfrak{sl}(2,\mathbb{R})\to\mathfrak{sl}(2,\mathbb{R})$, such that $A^*\left(\frac{1}{8}B\right)=g_1$ and $A^*R^{\frac{1}{8}B}=R^{g_1}$. This transformation is explicitly not required to be an element of the Lie algebra of $\mathrm{Aut}(\mathrm{SL}(2,\mathbb{R})\times\mathrm{SL}(2,\mathbb{R}))$. A simple ansatz for $A$ is
					\begin{center}
						\scalebox{0.75}{
							$A:=\left(
								\begin{array}{ccc|ccc}
									1 & 0 & 0 & -1 & 0 & 0 \\
									0 & 1 & 0 & 0 & 1 & 0 \\
									0 & 0 & 1 & 0 & 0 & 1 \\
									\hline
									\overset{}{\vphantom{|}}0 & 0 & 0 & 1 & 0 & 0 \\
									0 & 0 & 0 & 0 & 1 & 0 \\
									0 & 0 & 0 & 0 & 0 & 1
								\end{array}
							\right)$.
						}
					\end{center}
				With the help of our Maple script \href{http://www.math.uni-hamburg.de/home/cortes/ad_hoc_isom_check_g1_B.txt}{\texttt{ad\char`_hoc\char`_isom\char`_check\char`_g1\char`_B.txt}} we verify that the above transformation indeed fulfils the requirements.
			\end{proof}
		\end{Lem}

		\begin{Lem}\label{lem_inequivalence_metrics}
			The left-invariant metrics corresponding to $\frac{1}{8}B$ and $g_2$ are not isometric.
			\begin{proof}
				This follows from the fact that $\nabla^{g_2}R^{g_2}$ does not vanish identically. This can be checked with Maple using our script \href{http://www.math.uni-hamburg.de/home/cortes/nabla_R_g2.txt}{\texttt{nabla\char`_R\char`_g2.txt.}}
			\end{proof}
		\end{Lem}

		\begin{Lem}\label{lem_widetildeg_nK}
			The left-invariant metric corresponding to $g_2$ is the up to homothety unique nearly pseudo-K\"ahler metric on $\mathrm{SL}(2,\mathbb{R})\times\mathrm{SL}(2,\mathbb{R})$.
			\begin{proof}
				The existence and uniqueness of a nearly pseudo-K\"ahler metric on $\mathrm{SL}(2,\mathbb{R})\times\mathrm{SL}(2,\mathbb{R})$ up to homothety has been proven in \cite{SS}. Furthermore, it has stabilizer $\mathrm{SL}(2,\mathbb{R})$, meaning that it must be part of our classification of left-invariant pseudo-Riemannian Einstein metrics bi-invariant under at least a one-parameter subgroup. We have seen in the proof of Lemma \ref{lem_inequivalence_metrics} that $\nabla^{g_2}R^{g_2}$ does not vanish identically, and this is also true for the nearly pseudo-K\"ahler metric. For an explicit check of that statement we can use our Maple script \href{http://www.math.uni-hamburg.de/home/cortes/nabla_R_nK.txt}{\texttt{nabla\char`_R\char`_nK.txt}}. Hence the left-invariant metric corresponding to $g_2$ is homothetic to the nearly pseudo-K\"ahler metric as claimed.
			\end{proof}
		\end{Lem}
		
		Together, Lemmas \ref{lem_isometry_class_solutions}, \ref{lem_family_solutions_representative}, \ref{lem_equiv_B_g1}, \ref{lem_inequivalence_metrics}, and \ref{lem_widetildeg_nK} show that Theorem \ref{thm_main} holds true.
		\vspace{4em}
		
		\noindent
		\textbf{Conflicts of interest/Competing interests statement}\\
		The authors have no conflicts of interest to declare that are relevant to the content of this article.\\

		\noindent
		\textbf{Data availability statement}\\
		Some computations presented in the document have been carried out with the help of Maple  2021, using the license of the University of Hamburg. At
			\begin{center}
				\href{http://www.math.uni-hamburg.de/home/cortes/CEHL_Maple_files.zip}{\texttt{http://www.math.uni-hamburg.de/home/cortes/CEHL\char`_Maple\char`_files.zip}}
			\end{center}
		you can download the Maple files that we have used.

\end{document}